\documentclass[a4paper,11pt]{article}
\usepackage{amssymb}
\usepackage{latexsym}
\newtheorem{theorem}{Theorem}[subsection]

\newtheorem{proposition}{Proposition}[subsection]
\newtheorem{Lemma}{Lemma}[section]
\newtheorem{lemma}{Lemma}[subsection]

\newtheorem{remark}{Remark}[subsection]

\newcommand{\dis}{\displaystyle}

\begin{document}
\title{Asymptotic behavior of solutions to the Cauchy problem for 1-D $p$-system with space dependent damping}
\author{Akitaka Matsumura  and   Kenji Nishihara
\footnote{E-mail addresses: akitaka@muh.biglobe.ne.jp (A. Matsumura), kenji@waseda.jp (K. Nishihara)}
 \\[10pt]
{\small {\it Osaka University and Waseda University, Japan}}}
\date{}
\maketitle

\begin{abstract}
We consider the Cauchy problem for one-dimensional $p$-system with damping of space-dependent coefficient. This system models the  compressible flow through porous media in the Lagrangean coordinate. Our concern is an asymptotic behavior of solutions, which is expected to be the diffusion wave based on the Darcy law. In fact, in the constant coefficient case Hsiao and Liu \cite{hsiao-liu-92} showed the asymptotic behavior  under suitable smallness conditions for the first time. After this work, there are many literatures, but there are few works in the space-dependent damping case, as far as we know. In this paper we treat this space-dependent case, as a first step when the coefficient is around some positive constant.

\end{abstract}

\section{Introduction}
We consider the Cauchy problem for the $p$-system with space-dependent damping
$$
 \left\{ \begin{array}{l}
 v_t - u_x=0, \ \ (t,x) \in  {\mathbf R}_+ \times {\mathbf R}, \\[3pt] 
 u_t + p(v)_x = -\alpha u \ \ (\alpha = \alpha(x) ), \\[3pt]
 (v,u)(0,x) = (v_0,u_0)(x) \to  (v_{\pm},u_{\pm}), \ x \to \pm \infty, \ (v_{\pm}>0),
 \end{array}\right.
 \leqno{(1.1)} 
$$
which models the 1-D compressible flow through porous media in the Lagrangean coordinate, where $u=u(t,x)$ is the velocity of the flow at time $t$ and position $x$, $v\,(>0)$ is the specific volume, $p\,(>0)$ is the pressure with $p'(v)<0$, and the coefficient $\alpha = \alpha(x)>0$. Our interest is in the large time behavior of the solution $(v,u)$ to (1.1), which is expected to be the diffusion wave $(\bar{v},\bar{u})$ to
$$
 \left\{ \begin{array}{l}
   \bar{v}_t - \bar{u}_x = 0, \\[3pt]
   \quad p(\bar{v})_x = - \alpha \bar{u},
  \end{array}\right.
  \leqno{(1.2)}
 $$
 by the Darcy law. In fact, when $\alpha$ is a constant, Hsiao and Liu \cite{hsiao-liu-92} showed its asymptotic  behavior to $(\bar{v},\bar{u})$ under some smallness conditions for the first time. The convergence rates of $(v-\bar{v}, u-\bar{u})$ were improved by the second author \cite{nishihara-96}. After these, there are many works including the case $\alpha = \alpha(t)$. We cite only 
\cite{huang-marcati-pan-05, ji-mei-23, marcati-mei-rubino-05, mei-10, nishihara-wang-yang-00, zhang-22, zhang-zhu-23, zhao-01}. See the references therein, too. However, as far as we know, there are few results on the asymptotics toward diffusion wave in space-dependent damping case. Concerning the blow-up results, there are some works\,(see Chen et al. \cite{chen-li-li-mei-zhang-20}, Sugiyama \cite{sugiyama-23} and reference therein).

In this paper we consider the case of $x$-dependent coefficient $\alpha=\alpha(x)$. The Darcy law says that the velocity of flow in porous media is proportional to the pressure gradient, and so it may be reasonable that the coefficient $\alpha$ depends on the position $x$. By the situation of media, $\alpha$ will be assumed variously. However, as a first step we treat the simplest case
$$
  \alpha = \alpha(x) \to \underline{\alpha}, \ \ x \to \pm \infty
  \leqno{(1.3)}
$$
for some constant $\underline{\alpha}>0$. When 
$$
 v_+ = v_- \ \ (\mbox{and} \ u_+ \ne u_- \ \mbox{in general}),
 \leqno{(1.4)}
$$
the asymptotic profile $\bar{v}$ of (1.2) is expected to be a variant of the Gauss function. While, when
$$
  v_+ \ne v_-,  \leqno{(1.5)}
$$
$\bar{v}$ be an approximate similarity solution, and the treatment becomes complicated. 

In both cases, the problems are reformulated to the Cauchy problems for the quasilinear wave equations with damping of the coefficient $\alpha(x)$, whose results are directly applied to the original ones. Those details are stated in the next section. As for the results obtained in this paper we believe that there are still plenty of rooms for improvement. Under various assumptions on $\alpha = \alpha(x)$ or $\alpha(t,x)$ many discussions are expected in the future.
\\

The content of this paper is as follows. In Section 2 our problem is reformulated in each case of (1.4) and (1.5) to the Cauchy problem for a quasilinear wave equation of second order with damping. The case around the constant state of $v$ is treated in Subsection 2.1, while the different end states case is mentioned in Subsection 2.2. In each subsection, the theorem on the reformulated one will be stated clear, so that our goal for (1.1) in each case is obtained as a corollary. For those proofs, a standard energy method is applied, for which a series of a priori estimates are necessary. The case around the constant state (1.4) is shown in Section 3. In the final section the convergence to the similarity solution in case of (1.5) will be treated. \\

{\it Notations}. For the function spaces, $L^p= L^p({\mathbf R}) \, (1\le p \le \infty)$ is a usual Lebesgue space with norm
$$
 |f|_p = ( \int_{{\mathbf R}} |f(x)|^p \, dx)^{1/p} \ (1\le p < \infty), \quad |f|_{\infty} = \sup_{{\mathbf R}} |f(x)|.
$$
The integral domain ${\mathbf R}$ is often abbreviated when it is clear. For any integer $l \ge 0$, $H^l = H^l({\mathbf R})$ denotes the usual $l$-th order Sobolev space with norm
$$
   \| f\|_{H^l} = \| f \|_l = ( \sum_{j=0}^l \| \partial_x^j f \|^2)^{1/2}.
$$
When $l=0$ and $p=2$, we often use the notation $\| f \| = \| f\|_0 = |f|_2$. For brevity, $\| f,g, \cdots\|_{H^k\times H^l \times \cdots}^2$ $= \|f\|_{H^k}^2 + \|g\|_{H^l}^2 + \cdots = \|f\|_k^2 +\|g\|_l^2 + \cdots$. The set of $k$-times continuously differentiable functions in ${\mathbf R}$ with compact support is denoted by $C^k_0({\mathbf R})$. The space $C^k ([0,T];\, X)$ is a set of $k$-times differentiable functions on $[0,T]$ to the Hilbert space $X$. Also, by $C$ or $c$ we denote a generic positive constant independent of the data and time $t$, whose value may change in line to line.

\section{Reformulation of the problem and results}

In this section we heuristically reformulate the problem (1.1) to the quasi-linear wave equation with damping, and state the results on the reformulated one, which derive our main results on (1.1).

Rewrite (1.1) as
$$
 \left\{ \begin{array}{l}
  v_t - u_x = 0, \\[3pt]
  u_t + p(v)_x + \alpha u =0, \\[3pt]
  (v,u)(0,x) = (v_0,u_0)(x) \to (v_{\pm}, u_{\pm}), \ x \to \pm \infty,
 \end{array} \right.
 \leqno(2.1)
$$
under the assumption (1.3). By the Darcy law, the asymptotic profile $(\bar{v},\bar{u})$ is expected to be given by (1.2) as in the previous papers. However, in the case that $\alpha$ is depending on $x$ we meet difficulties to directly construct $(\bar{v},\bar{u})$ by
$$
 \bar{v}_t - \Big(\frac{-p(\bar{v})_x}{\alpha}\Big)_x = 0 \ \ \mbox{and} \ \ \bar{u} = \frac{-p(\bar{v})_x}{\alpha}.
 \leqno{(2.2)}
$$
In this paper, to avoid the difficulties, we shall introduce a simpler profile in each case of  (1.4) and (1.5). 

\subsection{Reformulation of the problem around the constant state}

Start to discuss the simpler case (1.4). Since both (1.3) and (1.4) are assumed, we here adopt an asymptotic profile $(V,U)$ by
$$
 \left\{ \begin{array}{l}
  V_t - U_x = 0, \\[5pt]
  \quad \ \ p'(\underline{v})V_x + \underline{\alpha} U = 0, \ V(t,\pm\infty) =\underline{v}\,( := v_+ = v_-),
  \end{array}\right.
  \leqno{(2.1.1)}
$$
or
$$
 \left\{ \begin{array}{l}
  V_t - \mu V_{xx} = 0, \ V(t,\pm\infty)= \underline{v} \ (\mu := \frac{|p'(\underline{v})|}{\underline{\alpha}}), \\[5pt]
 U= \frac{|p'(\underline{v})|}{\underline{\alpha}} V_x = \mu V_x,
 \end{array} \right.
 \leqno{(2.1.2)}
$$
so that, as an exact solution,
$$
 \left\{ \begin{array}{l}
  V(t,x) = \underline{v} + \frac{\delta_0}{\sqrt{4 \pi \mu (1+t)}} e^{- \frac{x^2}{4\mu (1+t)}}, \\[5pt]
  U(t,x) = \frac{-2\pi \delta_0 x}{(4\pi \mu (1+t))^{3/2}} e^{- \frac{x^2}{4\mu (1+t)}},
  \end{array} \right.
  \leqno{(2.1.3)}
$$
where $\delta_0$ is some constant, determined later. Then $(V,U)$ satisfies
$$
 \left\{ \begin{array}{l}
  V_t - U_x =0, \\[5pt]
  U_t + p(V)_x + \alpha U = U_t + (\alpha-\underline{\alpha}) U + (p(V) - p'(\underline{v})V)_x.
 \end{array} \right.
 \leqno{(2.1.4)}
$$
Here, we note that $(V,U) \to (\underline{v},0)$ as $x \to \pm \infty$.

To arrange the condition $u(0,x) = u_0(x) \to u_{\pm} \ (x\to \pm \infty)$, we introduce a correction function $(\hat{v},\hat{u})$. By (2.1)$_2$ (second equation of (2.1)) we can expect $u(t,x) \sim e^{-\underline{\alpha}t} u_{\pm}$ as $x\to \pm \infty$ and so we determine $(\hat{v},\hat{u})$ as
$$
 \left\{ \begin{array}{l}
  \hat{v}_t - \hat{u}_x = 0, \\[3pt]
  \hat{u}_t \quad\quad  + \alpha \hat{u} = 0, \\[3pt]
  (\hat{v},\hat{u})(t,x) \to (0, e^{-\underline{\alpha}t} u_{\pm}), \ x \to \pm \infty.
  \end{array} \right.
  \leqno{(2.1.5)}
$$
Thus we define $(\hat{v},\hat{u})$ by
$$
 \left\{ \begin{array}{l} 
  \hat{u}(t,x) = e^{-\alpha t} (u_- + (u_+ -u_-) \int_{-\infty}^x m_0(y)\, dy) =: e^{-\alpha t} M_0(x), \\[5pt]
  \hat{v}(t,x) = \frac{\partial}{\partial x} (\frac{e^{-\alpha t}}{-\alpha} M_0(x)) \\[5pt]
  \quad\quad\quad = \frac{e^{-\alpha t}}{-\alpha} \{ (u_+ -u_-) m_0(x) - \alpha'(x) M_0(x) (\frac{1}{\alpha(x)} + t) \},
  \end{array} \right.
  \leqno{(2.1.6)}
$$
where $m_0 \in C_0^{\infty}({\mathbf R})$ with $\int_{-\infty}^{\infty} m_0(y)\,dy = 1$. In the result, (2.1.5) holds with
$$
 \int_{-\infty}^{\infty} \hat{v}(0,x)\,dx = \int_{-\infty}^{\infty} \frac{\partial}{\partial x} (\frac{1}{-\alpha} M_0(x))\,dx = \frac{1}{-\underline{\alpha}} (u_+-u_-).
 \leqno{(2.1.7)}
$$

Combining (2.1), (2.1.4), (2.1.5), we have
$$
 \left\{ \begin{array}{l}
   (v-V-\hat{v})_t - (u-U-\hat{u})_x =0, \\[5pt]
   (u-U-\hat{u})_t + (p(v)-p(V))_x + \alpha (u-U-\hat{u}) \\[5pt]
   \quad = -\{ U_t + (\alpha -\underline{\alpha})U + (p(V)-p'(\underline{v})V)_x \}.
  \end{array} \right.
 \leqno{(2.1.8)}
$$
Since $u-U-\hat{u} \to 0 \ (x \to \pm \infty)$ is expected,
$$
 \begin{array}{l}
  \int_{-\infty}^{\infty} (v-V-\hat{v})(t,x)\,dx = \int_{-\infty}^{\infty} (v_0(x) -V(0,x) -\hat{v}(0,x))\,dx \\[5pt]
  \quad = \int_{-\infty}^{\infty} \{ (v_0(x) - \underline{v}) - (V(0,x) -\underline{v})\}\, dx + \frac{1}{\underline{\alpha}}(u_+ -u_-) \\[5pt]
  \quad = \int_{-\infty}^{\infty} (v_0(x) -\underline{v})\,dx - \delta_0 + \frac{1}{\underline{\alpha}}(u_+ -u_-) .
 \end{array}
 $$
Therefore, we choose $\delta_0$ by
$$ 
  \delta_0 = \int_{-\infty}^{\infty} (v_0(x) - \underline{v})\,dx + \frac{1}{\underline{\alpha}}(u_+ -u_-) ,
  \leqno{(2.1.9)}
$$
so that $\int_{-\infty}^{\infty} (v-V-\hat{v})(t,x)\,dx =0$ for any $t\ge 0$. Thus, if we define $\phi$ by
$$
  \phi(t,x) = \int_{-\infty}^{x} (v-V-\hat{v})(t,y)\,dy,
  \leqno{(2.1.10)}
$$
then we may expect $\phi(t,\cdot) \in H^1$, and it holds that $\phi_x = v-V-\hat{v}$. Then, by (2.1.8)$_1$, it also holds $\phi_t = u-U-\hat{u}$, and (2.1.8)$_2$ can be written as
$$
 \begin{array}{l}
   \phi_{tt} + (p(V+\hat{v}+\phi_x)-p(V))_x + \alpha \phi_t \\[5pt]
   \quad = - \{ U_t  + (\alpha-\underline{\alpha})U + (p(V)-p'(\underline{v})V)_x \}.
  \end{array}
  \leqno{(2.1.11)}
$$
Modify the second term in (2.1.11) to get the following reformulated problem:
$$
 \left\{ \begin{array}{l}
   \phi_{tt} + (p(V+\hat{v}+\phi_x)-p(V+\hat{v}))_x + \alpha \phi_t \\
  \quad = -\{ U_t + (\alpha-\underline{\alpha})U + (p(V)-p'(\underline{v})V)_x + (p(V+\hat{v})-p(V))_x \} =: F =: {\dis \sum_{i=1}^4 F_i}, \\[5pt]
  \phi(0,x) = \phi_0(x) := \int_{-\infty}^x (v_0(y) -V(0,y) -\hat{v}(0,y))\,dy, \\[5pt]
  \phi_t(0,x) = \phi_1(x) := u_0(x)-U(0,x) -\hat{u}(0,x).
 \end{array} \right.
 \leqno{(2.1.12)}
$$

Our aim is to show $(\phi_x,\phi_t)(t) = o(t^{-1/4}, t^{-3/4)}$ in $L^2$-sense and, if possible, $o(t^{-1/2}, t^{-1})$ in $L^{\infty}$-sense as $t\to \infty$, so that $(V,U+\hat{u})$ is an asymptotic profile of $(v,u)$, since $(\hat{v}, \hat{u})$ decays rapidly.

Thus, our next job is to show the existence and asymptotic behavior of the unique global-in-time solution $\phi$ to (2.1.12). We assume
$$
 \left\{ \begin{array}{l}
   p \in C^3({\mathbf R}_+), \ p'(v)<0 \ (v>0), \\[5pt]
   \alpha \in C^2({\mathbf R}), \ \alpha-\underline{\alpha} \in L^1\cap L^{2}, \ |x|^{1/2}(\alpha -\underline{\alpha}) \in L^2, \\[5pt]
   (1+|x|)(|\alpha_x|+|\alpha_{xx}|) \in L^2 \ \mbox{and} \ \alpha \ge \alpha_0>0\ (\alpha_0:\,\mbox{constant}),
 \end{array} \right.
 \leqno{(2.1.13)}
$$
and 
$$
 v_0 - \underline{v} \in L^1, \ \int_{-\infty}^x(v_0-\underline{v})(y)\,dy \in L^2, \ v_0-\underline{v} \in H^2 \ \mbox{and} \ u_0 -u_{\pm} \in L^2({\mathbf R}_{\pm}), u_{0x} \in H^1,
 \leqno{(2.1.14)}
$$
so $(\phi_0,\phi_1) \in H^3 \times H^2$. Under these assumptions there exists a unique local-in-time solution $\phi$ in $\cap_{i=0}^2 C^i([0,t_0];\,H^{3-i})$ for some $t_0>0$\,(for the proof, see Matsumura \cite{matsumura-77}). Main part in the proof for the global-in-time solution and its asymptotic behavior is the following a priori estimates.
\begin{proposition}[A priori estimate]
Assume (2.1.13) and (2.1.14). Then there exist positive constants $\varepsilon_0$ and $C_0$ such that, if $\delta_1:= |v_0-\underline{v}|_1+ |u_+-u_-| \le \varepsilon_0$ and if $\phi \in \cap_{i=0}^2 C^i([0,T];\,H^{3-i})$ is a solution to (2.1.12) for some $T>0$, which satisfies
$$
 \sup_{0\le t \le T} \{ \| (\phi, \phi_t)(t)\|_{H^3\times H^2} + (1+t) \| (\phi_{tx}, \phi_{txx})(t)\| \} \le \varepsilon_0,
$$
then it holds
$$
 \begin{array}{l}
  (1+t)^2 \| (\phi_{tt}, \phi_t, \phi_{xx})(t)\|_{H^1\times H^2\times H^1}^2 + (1+t) \| \phi_x(t)\|^2 + \| \phi(t)\|^2  \\[5pt]
  + \int_0^t \{ (1+\tau)^2 \| (\phi_{tt}, \phi_{tx}, \phi_{txx})(\tau)\|^2 + (1+\tau) \| (\phi_t,\phi_{xx})(\tau)\|^2 + \| \phi_x(\tau)\|^2 \}\, d\tau \\[5pt]
  \le C_0 (\| \phi_0,\phi_1 \|_{H^3 \times H^2}^2 + \delta_1), \ 0 \le t \le T.
 \end{array}
 \leqno{(2.1.15)}
$$
\end{proposition}

\begin{remark}
Note that none of smallness assumptions on $\alpha$ is not made.
\end{remark}

Combining the local existence and a priori estimates of solution, we have theorem for (2.1.12).

\begin{theorem}
Under the assumptions of Proposition 2.1.1, unique global-in-time solution $\phi \in  \cap_{i=0}^2 C^i([0,\infty);\, H^{3-i})$ to (2.1.12) exists and satisfies the decay properties(2.1.15) for $0\le t <\infty$.
\end{theorem}

Once we have Theorem 2.1.1, if we newly define $(v,u)$ by $(v,u) = (V+\hat{v}+\phi_x, U+\hat{u}+\phi_t)$, then we can have a unique solution $(v,u)$ to (2.1), satisfying $(v-V-\hat{v},u-U-\hat{u}) \in C^1([0,\infty);\,H^2)$. Thus, from Theorem 2.1.1 with 
$$
   |\phi_x(t)|_{\infty} \le C\| \phi_x(t) \|^{1/2} \| \phi_{xx}(t)\|^{1/2} \le C(1+t)^{-3/4},
$$
we have main theorem in case of (1.4).

\begin{theorem}[Around the constant state]
Suppose (2.1.13) and (2.1.14). If $\| \phi_0,\phi_1\|_{H^3 \times H^2}^2$ $+ \delta_1$ is suitably small, then the solution $(v,u)$ to (1.1) such that $(v-V-\hat{v},u -U-\hat{u}) \in C^1([0,\infty);\, H^2)$ uniquely exists and it holds that, as $t \to \infty$,
$$
\left\{ \begin{array}{l}
  (v- V)(t,x) = O(t^{-1/2}), \\[5pt]
  (u - U -\hat{u})(t,x) = O(t^{-1}),
 \end{array} \right.
 \quad \mbox{in} \quad L^2 \mbox{-sense},
 \leqno{(2.1.16)}
$$
and
$$
  v(t,x) = V(t,x) + O(t^{-3/4}) \quad \mbox{in} \quad L^{\infty} \mbox{-sense}.
  \leqno{(2.1.17)}
$$
\end{theorem}

Since $\| (V,U)(t)\|  \sim (t^{-1/4}, t^{-3/4})$ and $\|(\hat{v}, \hat{u})(t)\|_{H^1\times L^{\infty}} = O(e^{-ct})$, $(V,U+\hat{u})$ is an asymptotic profile of $(v,u)$ in $L^2$-framework, and $V$ is so of $v$ in $L^{\infty}$-framework.

\subsection{Reformulation around the different end states and results}

Similar to the last subsection, as an asymptotic profile we adopt the similarity solution $(V,U)$, which is defined by
$$
 \left\{ \begin{array}{l}
  V_t -U_x =0, \\[5pt]
  \quad p(V)_x + \underline{\alpha} U=0, \ \ V(t,\pm\infty)=v_{\pm} \ (v_+ \ne v_-),
  \end{array}\right.
  \leqno{(2.2.1)}
$$
or
$$
 \left\{ \begin{array}{l}
   V_t -(\mu(V)V_x)_x =0, \ \ V(t,\pm\infty)=v_{\pm}, \\[5pt]
   U= -\frac{1}{\underline{\alpha}} p(V)_x = -\frac{1}{\underline{\alpha}} p'(V) V_x =: \mu(V)V_x.
   \end{array}\right.
  \leqno{(2.2.2)}
$$
As a special solution, (2.2.2) has a similarity solution of the form
$$
  V(t,x) = \tilde{V}(\xi), \ \ \xi = (x-x_0)/\sqrt{1+t}.
  \leqno{(2.2.3)}
$$
where $x_0$ is any shift constant\,(cf. \cite{duyn-peletier-77}). We denote $\tilde{V}$ still by $V$.

\begin{lemma}
The similarity solution $V$ to (2.2.2)$_1$ satisfies the following:
$$
 \left\{ \begin{array}{l}
   |V(t,x)-v_{-}| \le C\delta_0 e^{-\frac{c(x-x_0)^2}{1+t}}, \ (t,x) \in {\mathbf R}_+ \times {\mathbf R}_-, \\[5pt]
   |V(t,x)-v_{+}| \le C\delta_0 e^{-\frac{c(x-x_0)^2}{1+t}}, \ (t,x) \in {\mathbf R}_+ \times {\mathbf R}_+,\\[5pt]
   |V_x(t,x)| \le C\delta_0 (1+t)^{-1/2} e^{-\frac{c(x-x_0)^2}{1+t}}, \ (t,x) \in {\mathbf R}_+ \times  {\mathbf R}, \\[5pt]
   |V_{xx}(t,x)| \le C\delta_0 (1+t)^{-1} e^{-\frac{c(x-x_0)^2}{1+t}}, \ (t,x) \in {\mathbf R}_+ \times  {\mathbf R}, \\[5pt]
   |V_{xxx}(t,x)| \le C\delta_0 (1+t)^{-3/2} e^{-\frac{c(x-x_0)^2}{1+t}}, \ (t,x) \in {\mathbf R}_+ \times {\mathbf R},
   \end{array}\right.
   \leqno{(2.2.4)}
$$
where $\delta_0 = |v_+ -v_-|$.
\end{lemma}

Since $V=V(\xi) \to v_{\pm}\, (\xi \to \pm \infty)$, $U=\mu(V)V_x \to 0$ as $x \to \pm\infty$. On the other hand, $u(0,x)=u_0(x) \to u_{\pm} \ (x \to \pm\infty)$, and so we must introduce the correction function $(\hat{v},\hat{u})$, samely as in the preceding subsection. We expect $u(t,x) \sim e^{-\underline{\alpha}t} u_{\pm}$ as $x \to \pm\infty$ by (2.1)$_2$, and so we determine $(\hat{v},\hat{u})$ as
$$
 \left\{ \begin{array}{l}
   \hat{v}_t -\hat{u}_x =0, \\[5pt]
   \hat{u}_t \quad +\alpha \hat{u}=0, \\[5pt]
   (\hat{v},\hat{u})(t,x) \to (0,e^{-\underline{\alpha}t} u_{\pm}), \ x \to \pm \infty.
   \end{array}\right.
   \leqno{(2.2.5)}
$$
Thus, we define $(\hat{v},\hat{u})$ by
$$
 \left\{ \begin{array}{l}
   \hat{u}(t,x) = e^{-\alpha t} (u_- + (u_+ -u_-)\int_{-\infty}^x m_0(y)\,dy) =: e^{-\alpha t}M_0(x), \\[5pt]
   \hat{v} = \frac{\partial}{\partial x} (\frac{e^{-\alpha t}}{-\alpha} M_0(x)) \\[5pt]
   \quad = \frac{e^{-\alpha t}}{-\alpha} \{ (u_+ -u_-) m_0(x) -\alpha'(x)M_0(x) (\frac{1}{\alpha(x)} + t)\},
   \end{array}\right.
   \leqno{(2.2.6)}
$$
where $m_0 \in C_0^{\infty} ({\mathbf R})$ with $\int_{-\infty}^{\infty} m_0(y)\,dy = 1$. By (2.1), (2.2.1) and (2.2.5),
$$
 \left\{ \begin{array}{l}
   (v-V-\hat{v})_t -(u-U-\hat{u})_x =0, \\[5pt]
   (u-U-\hat{u})_t + (p(v)-p(V))_x + \alpha (u-U-\hat{u}) = -\{ U_t + (\alpha -\underline{\alpha})U\}.
   \end{array}\right.
   \leqno{(2.2.7)}
$$
Since $u-U-\hat{u} \to 0$ is expected as $x\to\pm\infty$, integrating (2.2.7)$_1$ over $[0,t]\times {\mathbf R}$, we have
$$
 \int_{-\infty}^{\infty} (v-V-\hat{v})(t,x)\,dx = \int_{-\infty}^{\infty} (v_0(x)-V(x-x_0))\,dx + \frac{1}{\underline{\alpha}} (u_+ - u_-).
$$
For a given $v_0(x)$, $\int^{\infty}_{-\infty} (v_0(x) -V(x-x_0))\,dx \to \pm \infty$ or $\to \mp \infty$ ($x_0 \to \pm\infty$), and so we can choose $x_0$ as
$$
 \int_{-\infty}^{\infty} (v_0(x) -V(x-x_0))\,dx = -\frac{1}{\underline{\alpha}}(u_+-u_-),
 \leqno{(2.2.8)}
$$
(set $x_0=0$ without loss of generality) so that, for any $t \ge 0$,
$$
  \int_{-\infty}^{\infty} (v-V-\hat{v})(t,x)\,dx = 0.
$$
Thus, if we define $\phi$ by 
$$
 \phi(t,x) = \int_{-\infty}^x (v(t,y) -V(\frac{y}{\sqrt{1+t}}) -\hat{v}(t,y))\,dy,
 \leqno{(2.2.9)}
$$
then $(\phi_x,\phi_t) = (v-V-\hat{v}, u-U-\hat{u})$ and
$$
 \phi_{tt} +(p(V+\hat{v}+\phi_x)-p(V))_x +\alpha \phi_t = -\{ U_t + (\alpha -\underline{\alpha})U\}
$$
or, as in the previous subsection,
$$
 \begin{array}{l}
  \phi_{tt} + (p(V+\hat{v}+\phi_x)-p(V+\hat{v}))_x + \alpha \phi_t \\[5pt]
  \quad = -\{ U_t + (\alpha -\underline{\alpha})U + (p(V+\hat{v})-p(V))_x \} =: G =: {\dis \sum_{i=1}^3 G_i},
 \end{array} 
  \leqno{(2.2.10)}
$$
which is our reformulated problem in the present case. Formal reformulation is almost similar to the case around the constant state. However, since $V(\pm\infty)=v_{\pm}, \ v_+ \ne v_-$, the energy estimates of $\phi$ are rather different from those in case of $v_+ = v_-$. Our aim is to derive $\|(\phi_x,\phi_t)(t)\| = o(1,t^{-1/4})$ and, if possible, $|(\phi_x,\phi_t)(t)|_{\infty} = o(1, t^{-1/2})$ as $t\to\infty$, so that $(V,U+\hat{u})$ becomes an asymptotic profile of $(v,u)$ in $L^2$ and $L^{\infty}$ sense, respectively, since $(\hat{v},\hat{u})$ decays rapidly. Hence it is important to show the existence of global-in-time solution $\phi$ to (2.2.10), combining the existence of local-in-time solution and the a priori estimates with decay rates.

For the given data
$$
 (\phi,\phi_t)(0) = (\phi_0,\phi_1) \in H^3 \times H^2,
 \leqno{(2.2.11)}
$$
there exists a unique local-in-time solution $\phi \in \dis{\cap_{i=0}^2} C^i([0,t_0];\,H^{3-i})$ to (2.2.10) (for the proof, see Matsumura \cite{matsumura-77}). 

As usual, if we can have a priori estimates 
$$
  \| (\phi,\phi_t)(t)\|^2_{H^3 \times H^2} \le C(\| \phi_0,\phi_1\|^2_{H^3 \times H^2} + \delta_1), \ \delta_1 = |v_+ -v_-| + |u_+ - u_-|
$$
for the solution $\phi \in \cap_{i=0}^2 C^i([0,T];\,H^{3-i})$ with some decay properties, then continuation arguments are well done and global-in-time solution $\phi$ is obtained, which may satisfies desired decay rates. But, in this case it seems to be hopeless. In fact, to get the boundedness of $\| \phi(t)\|$, multiplying (2.2.10) by $\phi$ and integrating it over $[0,t] \times {\mathbf R}$, we have
$$
\int (\alpha \phi^2+\phi_t\phi)\,dx + \int_0^t \int (p(\hat{V})-p(\hat{V}+\phi_x))\phi_x \,dx\,d\tau
\le C(\|\phi_0,\phi_1\|^2 + \int_0^t \| \phi_t\|^2\, d\tau + \int_0^t \int G\phi\,dx\,d\tau),
$$
or
$$
 \| \phi(t)\|^2 + \int_0^t \|\phi_x(\tau)\|^2\,d\tau \le C(\|\phi_0,\phi_1\|^2 + \int_0^t \|\phi_t(\tau)\|^2\,d\tau + \int_0^t|\int G\phi\,dx|\,d\tau).
$$
where $\hat{V} = V+\hat{v}$. For an example, in the final term we need to estimate $G_2$ as
$$
 \int_0^t |\int (\alpha -\underline{\alpha})U\phi\,dx|\,d\tau \le \int_0^t |U|_{\infty} |\phi|_{\infty} |\alpha -\underline{\alpha}|_1\,d\tau \le C\delta_0\int_0^t (1+\tau)^{-1/2} \|\phi(\tau)\|^{1/2} \|\phi_x(\tau)\|^{1/2}\,d\tau,
$$
and to cover this by good terms, which seems to be impossible, even if the other terms can be well-evaluated.

However, fortunately there is not $\phi$ itself in the nonlinearity of (2.2.10), and so, for the
continuation arguments the uniform estimate in $[0,T]$ of 
$\| (\phi_x,\phi_t)(t)\|_2$ is necessary, while growing up of $\| \phi(t)\|$ may be permitted. 

Thus, our a priori estimate is the following.

\begin{proposition}[A priori estimate]
Assume that
$$
 \left\{ \begin{array}{l}
  p \in C^3({\mathbf R}_+), \ p'(v)<0 \, (v>0), \ \mbox{and} \ \alpha \in C^1({\mathbf R}), \ \alpha -\underline{\alpha} \in L^1 \cap L^2, \alpha_x \in L^2 \\[3pt]
  \mbox{with} \ \alpha \ge \alpha_0>0 \, (\alpha_0:\, \mbox{constant}).
 \end{array}  \right.
$$
Then there exist positive constants $\varepsilon_0$ and $C_0$ such that, if $\delta_1:= |v_+-v_-|+|u_+-u_-| \le \varepsilon_0$ and if $\phi \in \cap_{i=0}^2 C^i([0,T];\,H^{3-i})$ is a solution to (2.2.10)-(2.2.11) for some $T>0$, which satisfies
$$
 \sup_{0\le t\le T} \{ (1+t)^{-\gamma/2}\| \phi(t)\| + \| (\phi_t,\phi_x)(t)\|_2 + (1+t)^{1/2}\| (\phi_{txx}, \phi_{xxx})(t)\|\}  \le \varepsilon_0,
$$
then it holds that
$$
 \begin{array}{l}
   (1+t)^{-\gamma}\| \phi(t)\|^2 + (1+t)^{1-\gamma}\| \phi_x(t)\|^2 +(1+t)\| (\phi_t, \phi_{tt},\phi_{xx})(t)\|_{H^2\times H^1 \times H^1}^2 \\[5pt]

   \quad +\int_0^t \{ (1+\tau)^{-1-\gamma}\|\phi(\tau)\|^2 + (1+\tau)^{-\gamma}\|\phi_x(\tau)\|^2 + (1+\tau)^{1-\gamma}\|\phi_t(\tau)\|^2 \\[5pt]
   \quad \qquad + (1+\tau)\| (\phi_{tt},\phi_{tx})(\tau)\|^2 +\|(\phi_{ttx},\phi_{txx})(\tau)\|^2\}\,d\tau \\[5pt]
   \le C_0(\| \phi_0,\phi_1\|^2_{H^3 \times H^2} + \delta_1), \ 0\le t \le T.
 \end{array}
 \leqno{(2.2.12)}
$$
\end{proposition}

Proposition 2.2.1 yields the following theorem, together with the local existence theorem.

\begin{theorem}
Under the assumption in Proposition 2.2.1, if both $\| \phi_0,\phi_1\|_{H^3 \times H^2}$ and $\delta_1 = |v_+ -v_- | + |u_+-u_-|$ are sufficiently small, then there exists a unique global-in-time solution $\phi \in \cap_{i=0}^2 C^i([0,\infty);\,H^{3-i})$ to (2.2.10)-(2.2.11), which satisfies decay property (2.2.12) for $0\le t<\infty$.
\end{theorem}

As same as in the preceding subsection, once we have Theorem 2.2.1, if we define $(v,u)$ by $(v,u)=(V+\hat{v}+\phi_x,U+\hat{u}+\phi_t)$, then we can have a unique solution $(v,u)$ of (2.1) satisfying $(v-V-\hat{v},u-U-\hat{u}) \in C^1([0,\infty);\,H^2)$.

Also, since $\| \phi_x(t)\| \le C(1+t)^{-(1-\gamma)/2}$, $|\phi_x(t)|_{\infty} \le C\| \phi_x(t)\|^{1/2} \| \phi_{xx}(t)\|^{1/2} \le C(1+t)^{-(2-\gamma)/2}$ and $\| \phi_t(t)\| \le C(1+t)^{-1/2}$, the similarity solution $(V,U+\hat{u})$ is an asymptotic profile of $(v,u)$ in $L^2$-sense, and $V$ is so of $v$ in $L^{\infty}$-sense. Thus, our main theorem is the following.

\begin{theorem}[Main Theorem]
Under the assumption in Proposition 2.2.1, suppose that 
$$
  v_0-V \in L^1, \ \int_{-\infty}^{\cdot} (v_0(y)-V(0,y))\,dy \in L^2, \ (v_0 -V(0,\cdot),u_0-U(0,\cdot)-\hat{u}(0,\cdot)) \in H^2\times H^2.
$$
If both $\| \int_{-\infty}^{\cdot} (v_0(y)-V(0,y))\,dy\| + \| v_0 -V(0,\cdot), u_0-U(0,\cdot)-\hat{u}(0,\cdot)\|_{H^2 \times H^2}$ and $\delta_1 = |v_+ -v_-| + |u_+ -u_-|$ are suitably small, then there exists a unique solution $(v,u)$ to (1.1) satisfying $(v-V-\hat{v},u-U-\hat{u}) \in C^1([0,\infty);\,H^2)$, which satisfies, as $t\to\infty$, 
$$
 (v-V,u-U-\hat{u}) (t,x) = O(t^{-(1-\gamma)/2}, t^{-1/2}) \ \ \mbox{in} \ \ L^2\mbox{-sense},
$$
and
$$
  (v-V)(t,x) = O(t^{-(2-\gamma)/2}) \ \ \mbox{in} \ \ L^{\infty}\mbox{-sense}
$$
with $1/2 <\gamma<1$.
 \end{theorem}

\section{A priori estimate around the constant state}

We prove Proposition 2.1.1 in a series of several steps. To do that, let $\phi$ be a smooth solution in ${\dis \cap_{i=0}^2C^i([0,T];\,H^{3-i}})$ for some $T>0$ to (2.1.12) and rewrite (2.1.12) as
$$ 
  \phi_{tt} + (p(\hat{V}+\phi_x) -p(\hat{V}))_x + \alpha \phi_t = F,
  \leqno{(3.1)}
$$
where $\hat{V} := V+ \hat{v}$ and
$$
  F = -\{ U_t + (\alpha - \underline{\alpha})U + (p(V) -p'(\underline{v})V)_x + (p(V+\hat{v})-p(V))_x\}  =: \sum_{i=1}^4 F_i.
  \leqno{(3.2)}
$$
Also, denote the a priori assumption as
$$
 \delta = \sup_{0\le t \le T} \{ \| (\phi, \phi_t)(t)\|_{H^3\times H^2} + (1+t) \| (\phi_{tx}, \phi_{txx})(t)\| \} \,(\le 1) 
$$
together with
$$
  \delta_1 = |v_0-\underline{v}|_1+|u_+-u_-| \ge |\delta_0|+ (1+ \frac{1}{\underline{\alpha}})|u_+-u_-|
$$
and
$$
 I_0 = \| \phi_0,\phi_1\|_{H^3\times H^2}^2+\delta_1\,(\le 1),
$$
where $\delta_0$ is defined in (2.1.9). \\
\\
{\it Step 1.} First we multiply (3.1) by $\phi_t$ and integrate it over ${\mathbf R}$:
$$
  \begin{array}{l}
   \frac{d}{dt} \int \frac12 \phi_t^2 \,dx + \int \alpha \phi_t^2\,dx + \int (p(\hat{V})-p(\hat{V}+\phi_x))\phi_{tx} \,dx \\[5pt]
   \le \int |\phi_t F|\, dx \le \varepsilon \int \phi_t^2\,dx + C_{\varepsilon} \int F^2\, dx 
  \end{array}
  \leqno{(3.3)}
$$
for small constant $\varepsilon>0$. The integral domain ${\mathbf R}$ and variable $t$ are often abbreviated. The 3-rd term is estimated as
$$
 \begin{array}{l}
 \int (p(\hat{V})-p(\hat{V}+\phi_x)) \phi_{tx} \,dx \\[5pt]
 = \int (\int_0^{\phi_x} (p(\hat{V}) -p(\hat{V}+s))\,ds)_t\, dx + \int (p(\hat{V}+\phi_x) -p(\hat{V})- p'(V)\phi_x) \hat{V}_t\, dx \\[5pt]
 \ge \frac{d}{dt} \int (\int_0^{\phi_x} (p(\hat{V}) -p(\hat{V}+s))\,ds) \, dx - C(|V_t|_{\infty} + |\hat{v}_t|_{\infty} )\int \phi_x^2 \, dx \\[5pt]
  \ge \frac{d}{dt} \int (\int_0^{\phi_x} (p(\hat{V}) -p(\hat{V}+s))\,ds) \, dx - C\delta_1 \delta (1+t)^{-3/2} .
 \end{array}
$$

The estimates of $F$ are as follow, including those necessary. 

\begin{Lemma}[Estimates of $F$]
It holds that, for $t \in {\mathbf R}_+$
$$
 \left\{ \begin{array}{l}
   \|F(t)\|^2  \le C\delta_1(1+t)^{-5/2}, \\[5pt]
   \|F_x(t)\|^2  \le C\delta_1(1+t)^{-3}, \\[5pt]
   \|F_t(t)\|^2  \le C\delta_1(1+t)^{-9/2}, \\[5pt]
  \|F_{xt}(t)\|^2  \le C\delta_1(1+t)^{-5}.
\end{array} \right.
\leqno{(3.4)}
$$
\end{Lemma}

The proof is given in Appendix.

Thus, integrating (3.3) over $[0,t]$, we have 
$$
 \| (\phi_t,\phi_x)(t) \|^2 + \int_0^t \| \phi_t(\tau)\|^2 \, d\tau \le C \|\phi_{0x},\phi_1\|^2 + C(\delta_1 + \delta \delta_1)  
 \le CI_0.
 \leqno{(3.5)}
$$

Next, multiplying (3.1) by $\phi$, we have
$$
\begin{array}{l}
    \frac{d}{dt} \int (\frac{\alpha}{2}\phi^2 + \phi \phi_t) \,dx - \int \phi_t^2\, dx + \int (p(\hat{V})-p(\hat{V}+\phi_x))\phi_x\, dx \\[5pt]
    \le \int |\phi F|\,dx \le {\dis \sup_{0\le t \le T}}\| \phi\| \| F\| \le C\delta_1(1+t)^{-5/4},
 \end{array}
$$
and hence, using (3.5), 
$$
 \begin{array}{ll}
 \| \phi(t)\|^2 + \int_0^t \| \phi_x(\tau)\|^2\, d\tau & \le C(\| \phi_0,\phi_{0x},\phi_1\|^2 + \| \phi_t(t)\|^2 + \int_0^t \| \phi_t(\tau)\|^2\,d\tau )+C\delta_1 \\[5pt]
  & \le CI_0.
  \end{array}
 \leqno{(3.6)}
$$
Using (3.6), we can multiply (3.3) by $1+t$ and get
$$
 \| \phi(t)\|^2 + (1+t)\| (\phi_t,\phi_x)(t)\|^2 + \int_0^t (\| \phi_x(\tau)\|^2 + (1+\tau)\| \phi_t(\tau)\|^2)\,d\tau \le CI_0.
 \leqno{(3.7)}
 $$

We repeat the similar procedure to Step 1 for higher derivatives of $\phi$.
\\
\\
{\it Step 2.}  Differentiate (3.1) in $t$:
$$
 \phi_{ttt} + \alpha \phi_{tt} + (p'(\hat{V}+\phi_x)\phi_{xt} + (p'(\hat{V}+\phi_x)-p'(\hat{V}))\hat{V}_t)_x =F_t.
 \leqno{(3.8)}
$$
Multiplying (3.8) by $\phi_{tt}$ and integrating it over ${\mathbf R}$, we have
$$
 \begin{array}{l}
  \frac{d}{dt} \int ( \frac12 \phi_{tt}^2 + \frac{-p'(\hat{V}+\phi_x)}{2} \phi_{xt}^2)\,dx + \int \alpha \phi_{tt}^2 \,dx \\[5pt]
  \le C\delta_1 \{ (1+t)^{-4} \| \phi_x\|^2 + (1+t)^{-3} \| \phi_{xx}\|^2 + (1+t)^{-9/2}\} + C(\delta+\delta_1)\| \phi_{xt}\|^2,
  \end{array}
  \leqno{(3.9)}
$$
because of (3.4)$_3$, where the third term of (3.8) is estimated as
$$
 \begin{array}{l}
  \int (-p'(\hat{V}+\phi_x)\phi_{xt}\phi_{xtt} + \{ (p'(\hat{V}+\phi_x)-p'(\hat{V}))\hat{V}_{xt} + p''(\hat{V}+\phi_x)\hat{V}_t\phi_{xx} \\[5pt]
\quad + (p''(\hat{V} + \phi_x)-p''(\hat{V}))\hat{V}_x\hat{V}_t\} \phi_{tt})\, dx \\[5pt]
  \ge \frac{d}{dt} \int \frac{-p'(\hat{V}+\phi_x)}{2} \phi_{xt}^2\,dx + \int \frac{p''(\hat{V}+\phi_x)}{2} (\hat{V}_t+\phi_{xt})\phi_{xt}^2 \,dx \\[5pt]
  \quad -\varepsilon \int \phi_{tt}^2\,dx - C_{\varepsilon} \int (\hat{V}_{xt}^2 \phi_x^2 + \hat{V}_t^2 \phi_{xx}^2 +\hat{V}_x^2\hat{V}_t^2\phi_x^2)\,dx 
\end{array}
$$
$(0<\varepsilon \ll 1)$ and 
$$
 |\int \frac{p''(\hat{V}+\phi_x)}{2} (\hat{V}_t+\phi_{xt})\phi_{xt}^2 \,dx| \le C(\delta_1(1+t)^{-3/2}+\delta)\| \phi_{xt}\|^2  \le C(\delta_1+ \delta) \| \phi_{xt}\|^2
$$
by the a priori assumption. 

Also, the product of (3.8) and $\phi_t$ yields
$$
 \begin{array}{l}
  \frac{d}{dt} \int (\frac{\alpha}{2} \phi_t^2 + \phi_t \phi_{tt})\,dx + \frac{c}{2} \int \phi_{xt}^2\,dx \\[5pt]
  \le \int \phi_{tt}^2\,dx + C\delta_1 \{ (1+t)^{-9/4}\| \phi_t \|^2 + (1+t)^{-3}\| \phi_x \|^2 \},
 \end{array}
 \leqno{(3.10)}
$$
by $|\int \phi_t F_t \,dx| \le C\|\phi_t\| \|F_t\| \le C\delta_1 (1+t)^{-9/4}\| \phi_t\|$. By calculating (3.9)$+\lambda \cdot $(3.10)\, $(0<\lambda \ll 1)$,
$$
 \begin{array}{l}
  \frac{d}{dt} \int \{ (\frac12 \phi_{tt}^2 + \lambda \phi_{tt} \phi_t + \frac{\lambda \alpha}{2} \phi_t^2) + \frac{-p'(\hat{V}+\phi_x)}{2} \phi_{xt}^2 \}\,dx 
   + \int (\frac{\alpha}{2} \phi_{tt}^2 + \frac{\lambda c}{2} \phi_{tx}^2) \,dx \\[5pt]
  \le C(\delta +\delta_1 (1+t)^{-1}) \| \phi_{xt} \|^2 \\[5pt]
  \quad + C\delta_1 \{ (1+t)^{-3} \| \phi_x\|^2 + (1+t)^{-9/4} \| \phi_t\|^2 + (1+t)^{-9/2} + (1+t)^{-3} \| \phi_{xx}\|^2 \}.
 \end{array}
 \leqno{(3.11)}
$$
When $\delta + \delta_1 \ll 1$, integrating (3.11), $(1+t)\cdot$(3.11) and $(1+t)^2\cdot$(3.11) over $[0,t]$, yields
$$
 (1+t) \| (\phi_{tt},\phi_{tx},\phi_t)(t)\|^2 + \int_0^t (1+\tau) \| (\phi_{tt},\phi_{tx})(\tau)\|^2\,d\tau \le CI_0
 \leqno{(3.12)}
$$
and
$$
  (1+t)^2 \| (\phi_{tt},\phi_{tx},\phi_t)(t)\|^2 + \int_0^t (1+\tau)^2 \| (\phi_{tt},\phi_{tx})(\tau)\|^2\,d\tau \le CI_0 + C\delta_1 \int_0^t \| \phi_{xx}(\tau)\|^2\,d\tau.
  \leqno{(3.13)}
$$
\\
{\it Step 3.} In a similar fashion to $x$-derivative of (3.1), we have
$$
 \begin{array}{l}
   \frac{d}{dt} \int (\frac12 \phi_{tx}^2 + \frac{-p'(\hat{V}+\phi_x)}{2} \phi_{xx}^2) \,dx + c \int \phi_{tx}^2 \,dx \\[5pt]
   \le C\delta_1 \{ (1+t)^{-3} \| \phi_x \|^2 + (1+t)^{-3}\} + C\| \phi_t \|^2 + C\delta \|\phi_{xx}\|^2.
   \end{array}
   \leqno{(3.14)}
$$
Here, since we have differentiate (3.1) in $x$, the additional term $\alpha_x\phi_t$ comes out and
$$
 |\int \alpha_x \phi_t \phi_{tx}\,dx | \le \varepsilon \| \phi_{tx}\|^2 + C_{\varepsilon} \| \alpha_x\| \| \alpha_{xx}\| \cdot \| \phi_t\|^2.
$$
Though we omit the estimates of the other terms, integrating (3.14) and $(1+t)\cdot$(3.14) over $[0,t]$ and using (3.7), we have
$$
 (1+t) \| (\phi_{tx}, \phi_{xx})(t)\|^2 + \int_0^t (1+\tau)\| \phi_{tx}(\tau)\|^2\,d\tau \le CI_0 + C\delta \int_0^t \| \phi_{xx}(\tau)\|^2\,d\tau.
 \leqno{(3.15)}
$$

In (3.13) and (3.15), the final terms are not estimated yet.
\\
\\
{\it Step 4.} Multiplying the variant of (3.1)
$$
 \phi_{tt} + \alpha \phi_t + p'(\hat{V}+\phi_x)\phi_{xx} = F + (p'(\hat{V})-p'(\hat{V}+\phi_x))\hat{V}_x
 \leqno{(3.1)'}
$$
by $-\phi_{xx}$ and integrating it over ${\mathbf R}$, we obtain
$$
 \int \phi_{xx}^2 \,dx \le -\frac{d}{dt} \int \phi_x \phi_{xt}\,dx + \int \phi_{xt}^2\,dx +C\int \phi_t^2\,dx + C\delta_1 \{ (1+t)^{-2}\| \phi_x\|^2 + (1+t)^{-5/2}\}
$$
and, also integrating over $[0,t]$, 
$$
  \int_0^t \| \phi_{xx}(\tau)\|^2\,d\tau \le CI_0 + \| (\phi_t,\phi_{xt})(t)\|^2 + \int_0^t \|(\phi_t,\phi_{xt})(\tau)\|^2\,d\tau \le CI_0
  \leqno{(3.16)}
$$
by (3.12). Applying (3.16) to (3.13) and (3.15) with $\delta \ll 1$, we get
$$
 \begin{array}{l}
   (1+t)^2 \| (\phi_{tt},\phi_{tx},\phi_t)(t)\|^2 + (1+t)\| (\phi_x,\phi_{xx})(t)\|^2 +\| \phi(t)\|^2 \\[5pt]
   + \int_0^t  \{ (1+\tau)^2 \| (\phi_{tt},\phi_{tx})(\tau)\|^2 + (1+\tau)\| \phi_t(\tau)\|^2 + \|(\phi_{xx},\phi_x)(\tau)\|^2 \} \,d\tau \le CI_0.
   \end{array}
   \leqno{(3.17)}
$$
By (3.17), integration of $(1+t)^2(-\phi_{xx})\cdot$(3.1)$'$ in ${\mathbf R}$ and $(1+t)(-\phi_{xx})\cdot$(3.1)$'$ over $[0,t]\times {\mathbf R}$ also yield
$$
 (1+t)^2 \|\phi_{xx}(t)\|^2 + \int_0^t (1+\tau)\| \phi_{xx}(\tau)\|^2 \,d\tau \le CI_0.
 \leqno{(3.18)}
$$
\\
{\it Step 5.} Though the details are omitted, the combination of $\frac{\partial^2}{\partial t \partial x}$(3.1)$\cdot \phi_{ttx}$, $\frac{\partial}{\partial t}$(3.1)$\cdot (-\phi_{xxt})$ and $\frac{\partial}{\partial x}$(3.1)$\cdot (-\phi_{xxx})$ yields
$$
 (1+t)^2\|(\phi_{ttx}, \phi_{txx}, \phi_{xxx})(t)\|^2 + \int_0^t \{ (1+\tau)^2\| \phi_{txx}(\tau)\|^2 + (1+\tau)\|\phi_{xxx}(\tau)\|^2\} \,d\tau \le CI_0.
 \leqno{(3.19)}
$$

Thus, we have obtained (2.17) by (3.17) - (3.19) and completed the proof of Proposition 2.1.

\section{Proof of  Proposition 2.2.1 in different end states case}

Let $\phi$ be a smooth solution in $\cap_{i=0}^2 C^i([0,T];\,H^{3-i})$ for some $T>0$ to the reformulation problem (2.2.10)-(2.2.11), which is re-written as
$$
 \left\{ \begin{array}{l}
   \phi_{tt} + (p(\hat{V}+\phi_x)-p(\hat{V}))_x + \alpha \phi_t = G, \\[5pt]
   (\phi,\phi_t)(0) = (\phi_0,\phi_1) \in H^3 \times H^2,
 \end{array}\right.
 \leqno{(4.1)}
$$
where $\hat{V}=V + \hat{v}$ and
$$
 G= \sum_{i=1}^3 G_i = -\{ U_t + (\alpha-\underline{\alpha})U + (p(V+\hat{v})-p(V))_x \}.
 \leqno{(4.2)}
$$
To get the a priori estimate (2.2.12), denote the a priori assumption as
$$
 \delta = \sup_{0\le t\le T} \{ (1+t)^{-\gamma/2}\| \phi(t)\| + \| (\phi_t,\phi_x)(t)\|_2 + (1+t)^{1/2}\| (\phi_{txx}, \phi_{xxx})(t)\|\} \, (\le 1),
$$
together with
$$
 I_0 = \| \phi_0,\phi_1 \|_{H^3 \times H^2}^2 + \delta_1 \ (\delta_1= |v_+-v_-|+|u_+-u_-|).
$$
Note that the notations $\delta$ and $\delta_1$ in Sections 3 and 4 are slightly different from each other.

Same as in the preceding section, the proof is given in a series of several steps.
\\
\\
{\it Step 1}. Multiply (4.1) by $\phi_t$ and integrate it over ${\mathbf R}$:
$$
 \frac{d}{dt} \int \frac12 \phi_t^2 \,dx + \int \alpha \phi_t^2\,dx - \int (p(\hat{V}+\phi_x)-p(\hat{V}))\phi_{xt}\,dx = \int \phi_t G\,dx
 \leqno{(4.3)}
$$
and 
$$
 \begin{array}{l}
   \mbox{(3-rd term)} = \frac{d}{dt} \int \int_0^{\phi_x} (p(\hat{V})-p(\hat{V}+s))\,ds\,dx + \int (p(\hat{V}+\phi_x)-p(\hat{V})-p'(\hat{V})\phi_x)\hat{V}_t\,dx \\[5pt]
   \quad \ge \frac{d}{dt} \int \int_0^{\phi_x} (p(\hat{V})-p(\hat{V}+s))\,ds\,dx -C\delta_1(1+t)^{-1} \int \phi_x^2\,dx.
 \end{array}
$$
Hence, integrating (4.3) over $[0,t]$, we have
$$
  \begin{array}{l}
    \frac12 \| \phi_t(t)\|^2 + \int \int_0^{\phi_x} (p(\hat{V})-p(\hat{V}+s))\,ds\,dx + \int_0^t \int \alpha \phi_t^2\,dx\,d\tau \\[5pt]
    \le C\| \phi_0,\phi_1\|^2_{H^1 \times L^2} + C\delta_1 \int_0^t (1+\tau)^{-1}\|\phi_x(\tau)\|^2\,dx + \int_0^t\int \phi_t G\,dx\,dt,
  \end{array}
$$
or 
$$
 \|(\phi_t,\phi_x)(t)\|^2 + \int_0^t\|\phi_t(\tau)\|^2\,d\tau \le CI_0 + C\delta_1 \int_0^t (1+\tau)^{-1}\|\phi_x(\tau)\|^2\,dx + C\int_0^t\int \phi_t G\,dx\,dt.
 \leqno{(4.4)}
$$

Next, multiply (4.1) by $(1+t)^{-\gamma}\phi \, (1/2 <\gamma<1)$ and integrate it over $[0,t] \times {\mathbf R}$ to get
$$
 \begin{array}{l}
   (1+t)^{-\gamma}\int (\frac{\alpha}{2} \phi^2 + \phi \phi_t)\,dx + \gamma \int_0^t (1+\tau)^{-\gamma -1}\int (\phi^2 + \phi \phi_t)\,dx \,d\tau + c\int_0^t (1+\tau)^{-\gamma} \int \phi_x^2\,dx \,d\tau \\[5pt]
   \le CI_0 + \int_0^t (1+\tau)^{-\gamma}\int \phi_t^2 \,dx\,d\tau + C\int_0^t \int (1+\tau)^{-\gamma}\phi G\,dx\,d\tau.
 \end{array}
 \leqno{(4.5)}
$$
Because
$$
 \int \alpha (1+t)^{-\gamma}\phi \phi_t \,dx =\frac{d}{dt} (1+t)^{-\gamma}\int \frac{\alpha}{2} \phi^2 \,dx + \frac{\alpha \gamma}{2} (1+t)^{-\gamma-1} \int \phi^2\,dx, 
$$
$$
  \int (1+t)^{-\gamma} \phi \phi_{tt}\,dx = \frac{d}{dt} (1+t)^{-\gamma} \int \phi \phi_t\,dx + \gamma (1+t)^{-\gamma-1} \int \phi \phi_t \,dx - (1+t)^{-\gamma} \int \phi_t^2 \,dx,
$$
and
$$
 -\int (1+t)^{-\gamma} (p(\hat{V}+\phi_x)-p(\hat{V}))\phi_x \,dx \ge c(1+t)^{-\gamma}\int \phi_x^2\,dx.
$$

Adding (4.4) to $\lambda\cdot$(4.5) ($0< \lambda \ll 1$, in particular, $C\delta_1 \lambda \le c$), we have
$$
 \begin{array}{l}
  (1+t)^{-\gamma}\|\phi(t)\|^2 + \|(\phi_t,\phi_x)(t)\|^2  \\[5pt]
  \quad +\int_0^t ((1+\tau)^{-1-\gamma}\|\phi(\tau)\|^2 + (1+\tau)^{-\gamma}\| \phi_x(\tau)\|^2 + \|\phi_t(\tau)\|^2)\,d\tau \\[5pt]
  \le CI_0 + C \int_0^t (\int \phi_t G\,dx + (1+\tau)^{-\gamma} \int \phi G\,dx)\,d\tau.
 \end{array}
 \leqno{(4.6)}
$$
It is necessary to estimate the last term carefully. Estimate each term, using the a priori assumption (2.2.13) with $\delta_0=|v_+ -v_-|$:\\
\ \ (i) \ \ $\int \phi_t G_1\,dx \le C\delta_0 \int \phi_t^2\,dx + C\delta_0 \int (1+\tau)^{-3} e^{-\frac{2cx^2}{1+\tau}} \,dx \le C\delta_0 (\|\phi_t(\tau)\|^2 + (1+\tau)^{-5/2})$ \\[5pt]
by $|U_t| \le C(|V_x||V_t| + |V_{xt}|$,\\[10pt]
\ \ (ii) \ \ $\int \phi_tG_2\,dx = \frac{d}{dt} \int (\underline{\alpha}-\alpha)\phi U\,dx + \int (\alpha-\underline{\alpha})\phi U_t\,dx $ \\[5pt]
with
$$
  \begin{array}{l}
   |\int (\underline{\alpha} -\alpha)\phi U\,dx| \le C |U|_{\infty} \|\phi\| \| \alpha -\underline{\alpha}\| \le C\delta_0\delta (1+t)^{-(1-\gamma)/2} \le C\delta_0, \\[5pt]
   |\int (\alpha -\underline{\alpha} )\phi U_t\,dx|  \le |\phi|_{\infty}|U_t|_{\infty}|\alpha-\underline{\alpha}|_1 \le C\delta_0 (1+\tau)^{-3/2} \|\phi\|^{1/2} \|\phi_x\|^{1/2} \\[5pt]
   \quad \le C\delta_0 (1+\tau)^{(2\gamma-5)/4} \cdot (1+\tau)^{-(\gamma+1)/4}\|\phi\|^{1/2} \cdot (1+\tau)^{-\gamma/4} \|\phi_x\|^{1/2} \\[5pt]
   \quad \le C\delta_0\{ (1+ \tau)^{-(5-2\gamma)/2} + (1+\tau)^{-(\gamma+1)}\|\phi\|^2 + (1+\tau)^{-\gamma}\|\phi_x\|^2 \},
 \end{array}
$$
\ \ (iii) \ \ $\begin{array}{l}
              |(1+\tau)^{-\gamma}\int \phi G_1\,dx| \le C\delta_0 \int |\phi| (1+\tau)^{-\gamma-3/2} e^{-\frac{cx^2}{1+\tau}} \,dx \\
         \quad \le C\delta_0 ((1+\tau)^{-\gamma} \| \phi \|^2)^{1/2}               
(\int (1+\tau)^{-\gamma-3} e^{-\frac{2cx^2}{1+\tau}}\,dx)^{1/2} 
         \le C\delta_0\delta (1+\tau)^{-(2\gamma+5)/4} ,
    \end{array} $ \\[10pt]
\ \ (iv) \ \ $ \begin{array}{l}
               |(1+\tau)^{-\gamma} \int \phi (\underline{\alpha}-\alpha)U\,dx| \le (1+\tau)^{-\gamma} |\phi|_{\infty} |\alpha-\underline{\alpha}|_1 |U|_{\infty} \\[5pt]
              \quad \le C\delta_0 (1+\tau)^{-\gamma -1/2} \|\phi\|^{1/2} \|\phi_x\|^{1/2} \\[5pt]
              \quad =C\delta_0(1+\tau)^{-(2\gamma+1)/4} \cdot (1+\tau)^{-(\gamma+1)/4}\| \phi\|^{1/2} \cdot (1+\tau)^{-\gamma/4}\|\phi_x\|^{1/2} \\[5pt]
               \quad \le C\delta_0 \{ (1+\tau)^{-\gamma-1/2} + (1+\tau)^{-\gamma-1}\|\phi\|^2 + (1+\tau)^{-\gamma}\|\phi_x\|^2 \},
               \end{array} $\\[5pt]
and \\[5pt]
\ \ (v) \ \ since $|G_3| \le C|u_+-u_-| e^{-ct}$, $\int (\phi_t G_3 +(1+\tau)^{-\gamma}\phi G_3)\,dx $ is well estimated.\\

Thus, the final term in (4.6) is absorbed in the left-hand side if $\delta_1=\delta_0+|u_+-u_-|$ is small, and we have
$$
 \begin{array}{l}
  (1+t)^{-\gamma}\|\phi(t)\|^2 + \|(\phi_t,\phi_x)(t)\|^2  \\[5pt]
  \quad +\int_0^t ((1+\tau)^{-1-\gamma}\|\phi(\tau)\|^2 + (1+\tau)^{-\gamma}\| \phi_x(\tau)\|^2 + \|\phi_t(\tau)\|^2)\,d\tau 
  \le CI_0 .
 \end{array}
 \leqno{(4.7)}
$$
\\
{\it Step 2.}  We want to have decay properties of $\| (\phi_t,\phi_x)(t)\|$ and modify Step 1. In fact, multiplying (4.1) by $(1+t)^{1-\gamma}\phi_t$ and integrating it over $[0,t] \times {\mathbf R}$, we have
$$
 \begin{array}{l}
  (1+t)^{-\gamma}\|\phi(t)\|^2 + (1+t)^{1-\gamma}\|(\phi_t,\phi_x)(t)\|^2  \\[5pt]
  \quad +\int_0^t ((1+\tau)^{-1-\gamma}\|\phi(\tau)\|^2 + (1+\tau)^{-\gamma}\| \phi_x(\tau)\|^2 + (1+\tau)^{1-\gamma}\|\phi_t(\tau)\|^2)\,d\tau \\[5pt]
  \le CI_0,
 \end{array}
 \leqno{(4.8)}
$$
together with (4.6). Note that $1/2 <\gamma <1$, and so $(\phi_t,\phi_x)(t)$ decays in $L^2$-sense. 

Let us derive (4.8). Multiplying (4.3) by $(1+t)^{1-\gamma}$ yields
$$
 \begin{array}{l}
  \frac{d}{dt} \int \frac12 (1+t)^{1-\gamma}\phi_t^2 \,dx - \int \frac{1-\gamma}{2}(1+t)^{-\gamma} \phi_t^2\,dx + \int \alpha (1+t)^{1-\gamma}\phi_t^2\,dx \\[5pt]
  \quad +\frac{d}{dt} \int (1+t)^{1-\gamma} \int_0^{\phi_x} (p(\hat{V})-p(\hat{V}+s))\,ds\,dx - (1-\gamma) \int (1+t)^{-\gamma} \int_0^{\phi_x} (p(\hat{V})-p(\hat{V}+s))\,ds\,dx \\[5pt]
  \quad + (1+t)^{1-\gamma} \int (p(\hat{V}+\phi_x)-p(\hat{V}) -p'(\hat{V})\phi_x)\hat{V}_t\,dx \\[5pt]
= (1+t)^{1-\gamma} \int \phi_t G\,dx,
\end{array}
$$
and hence
$$
 \begin{array}{l}
   \frac{d}{dt} (1+t)^{1-\gamma} \{ \int \frac12 \phi_t^2 dx + \int \int_0^{\phi_x} (p(\hat{V})-p(\hat{V}+s))\,ds\,dx \} + \alpha_0 (1+t)^{1-\gamma} \int \phi_t^2\,dx \\[5pt]
   \le C(1+t)^{-\gamma} \|(\phi_t,\phi_x)(t)\|^2 + (1+t)^{1-\gamma}C\delta_1 (1+t)^{-1}\|\phi_x(t)\|^2 + C(1+t)^{1-\gamma} \int \phi_t G\,dx.
 \end{array}
$$
Here, the final term is estimated as follows and absorbed in the left-hand side:\\[5pt]
\ \ (vi) \ $(1+t)^{1-\gamma} \int \phi_t G_1\,dx \le (1+t)^{1-\gamma}\|\phi_t| \|U_t\| \le C\delta_0 \{ \|\phi_t(t)\|^2 +(1+t)^{-\gamma -3/2}\}$ \\[10pt]
\ \ (vii) \ $ \begin{array}{l}
                  (1+t)^{1-\gamma} \int \phi_t G_2\,dx 
                  = \frac{d}{dt} (1+t)^{1-\gamma} \int (\underline{\alpha}-\alpha)\phi U\,dx \\[5pt]
\quad- (1-\gamma)(1+t)^{-\gamma} \int (\underline{\alpha}-\alpha)\phi U\,dx -(1+t)^{1-\gamma} \int (\underline{\alpha}-\alpha)\phi U_t\,dx,
                 \end{array} $ \\[5pt]
 with \\[5pt]
 $ \begin{array}{l}
     (1+t)^{1-\gamma} \int (\underline{\alpha}-\alpha) \phi U\,dx \le (1+t)^{1-\gamma} |\phi|_{\infty} |U|_{\infty} |\alpha-\underline{\alpha}|_1 \\[5pt]
     \quad \le C\delta_0 (1+t)^{-\gamma+1/2} \| \phi\|^{1/2} \|\phi_x\|^{1/2} \\[5pt]
     \quad = C\delta_0 (1+t)^{-(2\gamma -1)/4} \cdot ((1+t)^{(1-\gamma)/2} \|\phi_x\|)^{1/2} \cdot ((1+t)^{-\gamma/2}\| \phi\|)^{1/2} \\[5pt]
     \quad \le C\delta_0 \{ (1+t)^{-\gamma} \|\phi(t)\|^2 + (1+t)^{1-\gamma}\| \phi_x(t)\|^2 + (1+t)^{-(2\gamma-1)/2}\}
     \end{array}
 $ \\[5pt]
 and the same estimates as (iii), (iv). \\[10pt]
 \ \ (viii)  For $G_3$, same as (v).\\[-5pt]
 
 Therefore, using (4.7), we have (4.8). \\
 
 We further need to evaluate higher derivatives of $\phi$.
 \\
 \\
{\it Step 3.}  Differentiate (4.1) in $t$:
$$
  \phi_{ttt} + \alpha \phi_{tt} + (p'(\hat{V}+\phi_x)\phi_{xt} + (p'(\hat{V}+\phi_x)-p'(\hat{V}))\hat{V}_t)_x = G_t. 
  \leqno{(4.9)}
$$
We multiply (4.9) by $\phi_t$, $\phi_{tt}$ and $(1+t)\phi_t$, $(1+t)\phi_{tt}$. First, estimate $\int \mbox{(4.9)}\cdot \phi_t\,dx$:
$$
 \begin{array}{l}
  \frac{d}{dt} \int (\phi_{tt}\phi_t + \frac{\alpha}{2} \phi_t^2)\,dx + \int |p'(\hat{V}+\phi_x)|\phi_{xt}^2\,dx \\[5pt]
  \quad -\int (p'(\hat{V}+\phi_x)-p'(\hat{V}))\hat{V}_t \phi_{tx}\,dx = \int \phi_t G_t\,dx.
  \end{array}
  \leqno{(4.10)}
$$
Since $|\hat{V}_t|_{\infty} \le C\delta_1 (1+t)^{-1} \le C\delta_1 (1+t)^{-\gamma}$,
$$
  |\mbox{3-rd term}| \le C\delta_1 \|\phi_{tx}(t)\|^2 + C\delta_1 (1+t)^{-\gamma}\| \phi_x(t)\|^2.
$$
It is easy to show
$$
  \|G_t\| \le C\delta_1 (1+t)^{-3/2}
  \leqno{(4.11)}
$$
and $|\int \phi_t G_t\,dx| \le C\delta_1(1+t)^{1-\gamma} \|\phi_t(t)\|^2 + C\delta_1(1+t)^{-(4-\gamma)}$, so that
$$
  \int_0^t \|\phi_{tx}(\tau)\|^2\,d\tau \le CI_0 + C(\|\phi_{tt}(t)\|^2 + \int_0^t\|\phi_{tt}(\tau)\|^2\,d\tau)
  \leqno{(4.12)}
$$
by (4.8). Secondly, multiplying (4.9) by $\phi_{tt}$,
$$
 \begin{array}{l}
   \frac{d}{dt} \int \frac12 \phi_{tt}^2\,dx + \int \alpha \phi_{tt}^2\,dx - \int p'(\hat{V}+\phi_x)\phi_{xt}\phi_{xtt}\,dx \\[5pt]
   \quad + \int ((p'(\hat{V}+\phi_x)-p'(\hat{V}))\hat{V}_t)_x\phi_{tt}\,dx = \int \phi_{tt} G_t\,dx,
 \end{array}
 \leqno{(4.13)}
$$
and 
$$
 \begin{array}{l}
  \mbox{(3-rd term)} =\frac{d}{dt} \int \frac{|p'(\hat{V}+\phi_x)|}{2} \phi_{xt}^2\,dx + \int \frac{p''(\hat{V}+\phi_x)}{2} (\hat{V}_t+\phi_{xt})\phi_{xt}^2\,dx \\[5pt]
  \quad \ge \frac{d}{dt} \int \frac{|p'(\hat{V}+\phi_x)|}{2} \phi_{xt}^2\,dx -C(\delta_1+\delta)\|\phi_{xt}(t)\|^2,
 \end{array}
$$
$$
 \begin{array}{l}
  |\mbox{4-th term}| \le C\delta_1 \|\phi_{tt}(t)\|^2 + C\delta_1 \int \{ (p'(\hat{V}+\phi_x)-p'(\hat{V}))\hat{V}_{xt} \\[5pt]
  \qquad + p''(\hat{V}+\phi_x)\phi_{xx}\hat{V}_t + (p''(\hat{V}+\phi_x)-p''(\hat{V}))\hat{V}_t\hat{V}_x\}^2\,dx \\[5pt]
  \quad \le C\delta_1\|\phi_{tt}(t)\|^2 + C\delta_1\{ (1+t)^{-3}\|\phi_x(t)\|^2 + (1+t)^{-2}\|\phi_{xx}(t)\|^2 \}. 
 \end{array}
$$
Here, by (4.1)
$$
  \phi_{xx} = (-p'(\hat{V}+\phi_x))^{-1} (\phi_{tt} +\alpha \phi_t + (p'(\hat{V}+\phi_x)-p'(\hat{V}))\hat{V}_x -G)
$$
and hence
$$
 (1+t)^{-2}\|\phi_{xx}(t)\|^2 \le C(1+t)^{-2} (\|(\phi_t,\phi_{tt})(t)\|^2 + \| G\|^2) + C\delta_1(1+t)^{-3}\|\phi_x(t)\|^2.
 \leqno{(4.14)}
$$
Since 
$$
 \|G\|^2 \le C(\|U_t\|^2 + |U|_{\infty}^2\|\alpha-\underline{\alpha}\|^2 + |u_+-u_-|e^{-ct}) \le C\delta_1 (1+t)^{-1}
  \leqno{(4.15)}
$$
and also
$$
 |\int \phi_{tt}G_t\,dx| \le \varepsilon \|\phi_{tt}(t)\|^2 + C_{\varepsilon}\int G_t^2\,dx \le  \varepsilon \|\phi_{tt}(t)\|^2 + C\delta_1(1+t)^{-3},
$$
all bad terms are absorbed, and integration of (4.13) over $[0,t]$ yields
$$
 \|(\phi_{tt},\phi_{tx})(t)\|^2 + \int_0^t \|\phi_{tt}(\tau)\|^2\,d\tau  \le CI_0 + C(\delta_1+ \delta)\int_0^t\| \phi_{xt}(\tau)\|^2\,d\tau.
 \leqno{(4.16)}
$$
Inserting (4.12) to (4.16), we have
$$
 \|(\phi_{tt},\phi_{tx})(t)\|^2 + \int_0^t \|(\phi_{tt},\phi_{tx}(\tau)\|^2\,d\tau \le CI_0,
 \leqno{(4.17)}
$$
provided that $\delta_1+\delta$ is suitably small.

Thirdly, based on (4.17), calculate $(1+t)\cdot$(4.10):
$$
  \begin{array}{l}
   \frac{d}{dt} (1+t) \int (\phi_{tt} \phi_t + \frac{\alpha}{2} \phi_t^2)\,dx + (1+t) \int |p'(\hat{V}+\phi_x)| \phi_{xt}^2\,dx \\[5pt]
   \le \int (\phi_{tt} \phi_t + \frac{\alpha}{2} \phi_t^2)\,dx + (1+t) \int (\phi_{tt}^2 + (p'(\hat{V}+\phi_x)-p'(\hat{V}))\hat{V}_t \phi_{tx})\,dx \\[5pt]
   \quad + (1+t) \int \phi_t G_t\,dx \\[5pt]
   \le C \int \phi_t^2\,dx + (1+t) \int (\phi_{tt}^2 + \varepsilon \phi_{tx}^2)\,dx \\[5pt]
   \quad + C_{\varepsilon} (1+t)|\hat{V}|_{\infty}^2 \int \phi_x^2\,dx + C(1+t)^{1-\gamma}\|\phi_t(t)\|^2 + C(1+t)^{1+\gamma}\|G_t\|^2   
   \end{array}
   \leqno{(4.18)}
$$
and
$$
 (1+t) \int \phi_{tt} \phi_t\,dx \le \varepsilon (1+t) \int \phi_t^2\,dx + C_{\varepsilon} (1+t) \int \phi_{tt}^2\,dx.
$$
Hence, integrating (4.18) over $[0,t]$ and using (4.8), (4.11), we have
$$
 \begin{array}{l}
   (1+t)\|\phi_t(t)\|^2 + \int_0^t (1+\tau) \|\phi_{xt}(\tau)\|^2\,d\tau \\[5pt]
   \le CI_0 + C(1+t) \|\phi_{tt}(t)\|^2 + C\int_0^t (1+\tau)\|\phi_{tt}(\tau)\|^2\,d\tau.
 \end{array}
 \leqno{(4.19)}
$$

Finally, multiplication of (4.13) by $1+t$ and integration over $[0,t]$ yield
$$
 \begin{array}{l}
  (1+t)\|(\phi_{tt},\phi_{tx})(t)\|^2 + \int_0^t (1+\tau)\|\phi_{tt}(\tau)\|^2\,d\tau \\[5pt]
  \le C \|(\phi_{tt},\phi_{tx})(t)\|^2 + C(\delta_1+\delta) \int_0^t (1+\tau)\|\phi_{tx}(\tau)\|^2 \,d\tau \\[5pt]
  \quad + C\delta_1 \int_0^t ((1+\tau) \|\phi_{tt}(\tau)\|^2 + (1+\tau)^{-2}\|\phi_x(\tau)\|^2 + (1+\tau)^{-1}\|\phi_{xx}(\tau)\|^2)\,d\tau.
 \end{array}
$$
By (4.14),
$$
  \|\phi_{xx}(t)\|^2 \le C\|(\phi_t,\phi_{tt})(t)\|^2 + C\delta_1(1+t)^{-1}\|\phi_x(t)\|^2 + C\delta_1(1+t)^{-1}.
$$
Therefore, by (4.17) we have
$$
 \begin{array}{l}
   (1+t) \|(\phi_{tt},\phi_{tx})(t)\|^2 + \int_0^t (1+\tau)\|\phi_{tt}(\tau)\|^2\,d\tau \\[5pt]
   \le CI_0 + C(\delta_1+\delta) \int_0^t (1+\tau)\|\phi_{tx}(\tau)\|^2\,d\tau.
 \end{array}
 \leqno{(4.20)}
$$
Substituting (4.20) into (4.19), we obtain
$$
 (1+t)\|(\phi_t,\phi_{tt},\phi_{tx})(t)\|^2 + \int_0^t (1+\tau) \|(\phi_{tt},\phi_{tx})(\tau)\|^2 \,d\tau \le CI_0,
 \leqno{(4.21)}
$$
provided that $\delta_1+\delta$ is small.

Additionally, we have
$$
 (1+t)\|\phi_{xx}(t)\|^2 \le CI_0,
 \leqno{(4.22)}
$$
because, by (4.14) and (4.15) with (4.21),
$$
 (1+t)\|\phi_{xx}(t)\|^2 \le C(1+t)\|(\phi_t,\phi_{tt})(t)\|^2 + C\delta_1\|\phi_x(t)\|^2 + C\delta_1.
$$
{\it Step 4.} In the final step we estimate the third order derivatives of $\phi$. To do so, differentiate (4.1) in $x$ and $x,t$:
$$
 \begin{array}{l}
  \phi_{ttx} + \alpha_x \phi_t + \alpha \phi_{tx} + (p'(\hat{V}+\phi_x)\phi_{xx})_x \\[5pt]
  = -((p'(\hat{V}+\phi_x)-p'(\hat{V}))_x + G_x,
  \end{array}
  \leqno{(4.23)}
$$
and
$$
 \begin{array}{l}
  \phi_{tttx} + \alpha_x \phi_{tt} + \alpha \phi_{ttx} + (p'(\hat{V}+\phi_x)\phi_{txx})_x \\[5pt]
  = -(p''(\hat{V}+\phi_x)(\hat{V}_t + \phi_{tx})\phi_{xx})_x -((p'(\hat{V}+\phi_x))\hat{V}_x)_{tx} + G_{tx} \\[5pt]
  =: H + G_{tx} =: h_1 + h_2 + G_{tx}.
  \end{array}
  \leqno{(4.24)}
$$

First, multiplying (4.24) by $\phi_{ttx}$ and integrating it over ${\mathbf R}$, we have
$$
 \begin{array}{l}
  \frac{d}{dt} \int (\frac{1}{2} \phi_{ttx}^2 + \frac{|p'(\hat{V}+\phi_x)|}{2} \phi_{txx}^2)\,dx + \int \alpha \phi_{ttx}^2\,dx \\[5pt]
  \le \varepsilon \int \phi_{ttx}^2\,dx + C_{\varepsilon} \int \phi_{tt}^2\,dx + C|\hat{V}_t+\phi_{tx}|_{\infty} \int \phi_{txx}^2\,dx + C_{\varepsilon} \int (H^2 + G_{tx}^2) \,dx,
 \end{array}
$$
which derives
$$
 \begin{array}{l}
  \frac{d}{dt} \int (\frac12 \phi_{ttx}^2 + \frac{|p'(\hat{V}+\phi_x)|}{2} \phi_{txx}^2)\,dx + \alpha_0 \int \phi_{ttx}^2\,dx \\[5pt]
  \le C\{ (\delta_1+\delta)\int \phi_{txx}^2\,dx + \int \phi_{tt}^2\,dx + \| H\|^2 + \| G_{tx}\|^2 \},
  \end{array}
  \leqno{(4.25)}
$$
since
$$
  |\hat{V}_t+\phi_{tx}|_{\infty} \le C(\delta_1+ \sqrt{I_0\delta}) (1+t)^{-1} \le C(\delta_1+\delta)(1+t)^{-1}
  \leqno{(4.26)}
$$
by (2.2.13), (4.21) and $\varepsilon \ll 1$. Next, multiplying (4.24) by $\phi_{tx}$, similarly as above, we have
$$
 \begin{array}{l}
  \frac{d}{dt} \int (\phi_{ttx}\phi_{tx} + \frac{\alpha}{2} \phi_{tx}^2)\,dx  -\int \phi_{ttx}^2 \,dx + \int |p'(\hat{V}+\phi_x)| \phi_{txx}^2 \,dx \\[5pt]
  \le C|\alpha_x|_{\infty}(\| \phi_{tt}\|^2 + \|\phi_{tx}\|^2) + C(\|\phi_{tx}\|^2 + \| H\|^2 + \| G_{tx}\|^2 ) \\[5pt]
  \le C\{ \| (\phi_{tt},\phi_{tx})(t)\|^2 + \| H\|^2 + \| G_{tx}\|^2\}.
  \end{array}
  \leqno{(4.27)}
$$
For small $\lambda>0$, add (4.25) to $\lambda\cdot$(4.27), and then
$$
 \begin{array}{l}
  \frac{d}{dt} \int \{ (\frac12 \phi_{ttx}^2 +\lambda \phi_{ttx}\phi_{tx} + \frac{\alpha\lambda}{2}\phi_{tx}^2) + \frac{|p'(\hat{V}+\phi_x)|}{2}\phi_{txx}^2\}\,dx \\[5pt]
  \quad + \int ( \frac{\alpha_0}{2}\phi_{ttx}^2 + \frac{\lambda}{2}|p'(\hat{V}+\phi_x)|\phi_{txx}^2)\,dx \\[5pt]
  \le C(\|(\phi_{tt},\phi_{tx})(t)\|^2 + \| H\|^2 + \| G_{tx}\|^2)
  \end{array}
  \leqno{(4.28)}
$$
provided that $\delta_1+\delta \le \frac{\lambda}{2}|p'(\hat{V}+\phi_x)|$. Here, $\| H\|^2$ is estimated as the following.
\begin{Lemma}
$$
 \begin{array}{l}
  \| H\|^2 \le C\{ (\delta_1+\delta) (1+t)^{-1} \| (\phi_{ttx},\phi_{txx},\phi_{tx},\phi_t)(t)\|^2 \\
  \qquad\qquad + (\delta_1+\delta)(1+t)^{-2}\| \phi_{xx}(t)\|^2 + \delta_1(1+t)^{-3}\|\phi_x(t)\|^2 + \delta_1 (1+t)^{-2} \}.
  \end{array}
  \leqno{(4.29)}
$$
\end{Lemma}
{\it Proof of Lemma 4.1}.  First, in (4.24)
$$
 \begin{array}{l}
  h_1 = -p''(\hat{V}+\phi_x)(\hat{V}_t +\phi_{tx})\phi_{xxx} - p''(\hat{V}+\phi_x)(\hat{V}_{tx}+\phi_{txx})\phi_{xx} \\[5pt]
  \qquad -p'''(\hat{V}+\phi_x)(\hat{V}_x+\phi_{xx})(\hat{V}_t+\phi_{tx})\phi_{xx} \\[5pt]
  =: h_{11}+h_{12}+h_{13},
 \end{array}
$$
with
$$
 \begin{array}{ll}
  \| h_{11}\|^2 & \le C(\delta_1+\delta)(1+t)^{-1} \| \phi_{xxx}\|^2 \ \ \mbox{by (4.26)} \\[5pt]
  \| h_{12}\|^2 & \le C(|\hat{V}_{tx}|_{\infty}^2 \| \phi_{xx}\|^2 + \|\phi_{xx}\| \|\phi_{xxx}\| \|\phi_{txx}\|^2) \\[5pt]
  & \le C\delta_1 (1+t)^{-3}\| \phi_{xx}\|^2 + (I_0+\delta)(1+t)^{-1}\| \phi_{txx}\|^2, \\[5pt]
  \|h_{13}\|^2 & \le C(|\hat{V}_x|_{\infty}^2 + \| \phi_{xx}\| \|\phi_{xxx}\|)(|\hat{V}_t|_{\infty}^2 + \|\phi_{tx}\| \|\phi_{txx}\|)\| \phi_{xx}\|^2 \\[5pt]
  & \le C(\delta_1 +\delta)(1+t)^{-2}\|\phi_{xx}\|^2.
 \end{array}
$$
For the estimate of $\|\phi_{xxx}\|^2$ in $\|h_{11}\|^2$ we back to (4.23):
$$
 \begin{array}{l}
   -p'(\hat{V}+\phi_x)\phi_{xxx} = p''(\hat{V}+\phi_x)(\hat{V}_x+\phi_{xx})\phi_{xx} + ((p'(\hat{V}+\phi_x)-p'(\hat{V}))\hat{V}_x)_x \\[5pt]
   \qquad\qquad + \phi_{ttx} + \alpha_x\phi_t +\alpha \phi_{tx} -G_x,
 \end{array}
$$
and hence
$$
 \begin{array}{ll}
  \|\phi_{xxx}\|^2 & \le C\{ (|\hat{V}_x|_{\infty}^2 + \|\phi_{xx}\| \|\phi_{xxx}\|)\|\phi_{xx}\|^2 + (|\hat{V}_{xx}|_{\infty}^2 + |\hat{V}_x|_{\infty}^4)\| \phi_x\|^2 \\[5pt]
   & \quad +|\hat{V}_x|_{\infty}^2 \|\phi_{xx}\|^2 + \| \phi_{ttx},\phi_{tx},\phi_t \|^2 + \|G_x\|^2\} \\[5pt]
   & \le C \{ (I_0+\delta)(1+t)^{-1}\|\phi_{xx}\|^2 + \delta_1(1+t)^{-2}\|\phi_x\|^2 \\[5pt]
  & \quad   + \| \phi_{ttx},\phi_{tx},\phi_t \|^2+\delta_1(1+t)^{-1}\}.
   \end{array}
   \leqno{(4.30)}
$$
Here, $\| G_x\|^2 \le C\delta_1 (1+t)^{-1}$ is easily seen. In a similar fashion to the above, we have
$$
 \|h_2\|^2 \le C\delta_1 \{ (1+t)^{-1}\|\phi_{txx}\|^2 + (1+t)^{-2}\|\phi_{tx},\phi_{xx}\|^2 +(1+t)^{-4}\|\phi_x\|^2\}.
$$
Combining $\|h_1\|^2, \|h_2\|^2$ with (4.30), we get (4.29) and complete the proof of Lemma 4.1.  \hfill {\it q.e.d.}
\\

We also note that
$$
 \| G_{tx}\|^2 \le C\delta_1 (1+t)^{-3}. 
 \leqno{(4.31)}
$$

We now return back to (4.28). Take $\delta_1+ \delta$ as sufficiently small, then the term $\phi_{txx}$ and $\phi_{ttx}$ in (4.29) are absorbed into the left hand side, and integration of (4.28) over $[0,t]$ yields
$$
 \| (\phi_{ttx},\phi_{txx},\phi_{tx})(t)\|^2 + \int_0^t \|(\phi_{ttx},\phi_{txx})(\tau)\|^2\,d\tau \le CI_0,
 \leqno{(4.32)}
$$
because of (4.8), (4.21)-(4.22).

To get further decay rate, we want to multiply (4.28) by $1+t$, but, if so, the final term in (4.29) is not integrable in $t$. So, we here use the technique found in Nishikawa \cite{nishikawa-98}, that is, multiply (4.28) by $(1+t)^{1+\nu} \, (\nu>0,$ not $\nu<0$), so that, by (4.32),
$$
 \begin{array}{l}
  (1+t)^{1+\nu} \| (\phi_{ttx},\phi_{txx},\phi_{tx})(t)\|^2 + \int_0^t (1+\tau)^{1+\nu}\|(\phi_{ttx},\phi_{txx})(\tau)\|^2\,d\tau \\[5pt]
  \le C\{ (1+\nu) \int_0^t(1+\tau)^{\nu}\|(\phi_{ttx},\phi_{txx},\phi_{tx})(\tau)\|^2\,d\tau + \int_0^t (1+\tau)^{1+\nu}\|(\phi_{tt},\phi_{tx})(\tau)\|^2\,d\tau \\[5pt]
  \quad + \int_0^t (1+\tau)^{1+\nu} \|H\|^2\,d\tau \} + \frac{C\delta_1}{\nu}(1+t)^{\nu}.
 \end{array}
 \leqno{(4.33)}
$$
Divide (4.33) by $(1+t)^{\nu}$, and use $\frac{1+\tau}{1+t} \le 1$ and (4.32) just obtained, then we get desired estimate
$$
 \begin{array}{l}
  (1+t)\| (\phi_{ttx},\phi_{txx},\phi_{tx})(t)\|^2 + \int_0^t \|(\phi_{ttx},\phi_{txx})(\tau)\|^2\,d\tau \\[5pt]
  \le C\{ \int_0^t \|(\phi_{ttx},\phi_{txx},\phi_{tx})(\tau)\|^2\,d\tau + \int_0^t (1+\tau) \|(\phi_{tt},\phi_{tx})(\tau)\|^2\,d\tau + \delta_1\} \\[5pt]
  \le CI_0.
  \end{array}
  \leqno{(4.34)}
$$
Here, note that, although $\frac{(1+\tau)^{1+\nu}}{(1+t)^{\nu}}$ comes out in the second term in (4.34), it only holds that $\frac{(1+\tau)^{1+\nu}}{(1+t)^{\nu}} \ge 0 \, (0\le \tau \le t)$.

Additionally, multiplying (4.30) by $1+t$, we get
$$
 (1+t)\| \phi_{xxx}(t)\|^2 \le CI_0.
 \leqno{(4.35)}
$$

Thus, we have obtained the estimate (2.2.12) and completed the proof of Proposition 2.2.1.
\\
\\
\\
{\Large\bf Appendix.} We prove Lemma 3.1. By
$$
 V-\underline{v} = \frac{\delta_0}{\sqrt{4\pi \mu (1+t)}} e^{-\frac{x^2}{4\mu (1+t)}},
 \leqno{(A1)}
$$
$$
 U=V_x = \frac{-2\pi \delta_0 x}{(4\pi \mu (1+t))^{3/2}} e^{-\frac{x^2}{4\mu (1+t)}},
 \leqno{(A2)}
$$
and (2.1.6), we easily know
$$
 \begin{array}{l}
  |\partial_x^k(V-\underline{v})|_p \le C \delta_0 (1+t)^{-\frac12(1-\frac{1}{p})-\frac{k}{2}}, \\[5pt]
  |\partial_x^k U|_p = |\partial_x^{k}V_x|_p \le C\delta_0 (1+t)^{-\frac12 (1-\frac{1}{p}) -\frac{k+1}{2}}
 \end{array}
 \leqno{(A3)}
$$
and
$$
 |\partial_x^k \hat{v}|_p + |\partial_t^l \hat{v}|_p \le C|u_+-u_-| e^{-ct}\le C\delta_1 e^{-ct} \ (0<c<\alpha_0).
 \leqno{(A4)}
$$
Rewrite (3.2):
$$
 F= -\{ U_t + (\alpha-\underline{\alpha})U +(p(V)-p'(\underline{v})V)_x + (p(V+\hat{v})-p(V))_x\} =: \sum_{i=1}^4 F_i.
 \leqno{(3.2)}
$$

By (A3)-(A4),
$$
 \begin{array}{l}
  \| F_1(t)\|^2 \le \| U_t\|^2 = \| V_{xxx}\|^2 \le C\delta_0 (1+t)^{-7/2}, \\[5pt]
  \| F_3(t)\|^2 \le \int |p'(V)-p'(\underline{v})|^2 |V_x|^2\,dx \le C|V_x|_{\infty}^2 \| V-\underline{v}\|^2 \le C\delta_0 (1+t)^{-5/2}, \\[5pt]
  \| F_4(t)\|^2 \le C\delta_1 e^{-ct}.
 \end{array}
$$
For $F_2$, by (A2)
$$
 \begin{array}{ll}
  \| F_2(t)\|^2 & \le \int |\alpha-\underline{\alpha}|^2|V_x|^2\,dx \le C\delta_0 \int \frac{|\alpha-\underline{\alpha}|^2 |x|}{(1+t)^{3/2}}\,dx \cdot |V_x|_{\infty} \\[5pt]
   & \le C\delta_0 (1+t)^{-5/2} \ \ \mbox{if} \ \ |x|^{1/2}(\alpha-\underline{\alpha}) \in L^2.
 \end{array}
$$
Hence we have $\| F(t)\|^2 \le C\delta_1(1+t)^{-5/2}$. Next,
$$
 \begin{array}{ll}
   \| F_{1x}(t)\|^2 & = \| U_{tx}\|^2 \le C\delta_0 (1+t)^{-9/2}, \\[5pt]
   \| F_{2x}(t) \|^2 & \le C \int (|\alpha_x U|^2 + |\alpha -\underline{\alpha}|^2 U_x^2)\,dx \\[5pt]
        & \le C\delta_0 (1+t)^{-3} (\| x\cdot \alpha_x\|^2 + \| \alpha-\underline{\alpha}\|^2) \le C\delta_0(1+t)^{-3}, \\[5pt]
   \| F_{3x}(t)\|^2 & \le \int (p''(V)V_x^2 + (p'(V)-p'(\underline{v}))V_{xx})^2\,dx \le C\delta_0(1+t)^{-7/2}, \\[5pt]
   \| F_{4x}(t)\|^2 & \le C\delta_1 e^{-ct},
 \end{array}
$$
and hence $\| F_x(t)\|^2 \le C\delta_1 (1+t)^{-3}$.  For $F_t$,
$$
 \begin{array}{ll}
  \| F_{1t}(t)\|^2 & = \| U_{tt}\|^2 = \| V_{xtt}\|^2 \le C\delta_0(1+t)^{-11/2}, \\[5pt]
  \| F_{3t}(t)\|^2 & = \| (p'(V)-p'(\underline{v}))V_{tx} + p''(V) V_xV_t\|^2 \\[5pt]
       & \le C\delta_0 ((1+t)^{-1/2 -4}+(1+t)^{-3/2 -3}) = C\delta_0 (1+t)^{-9/2}, \\[5pt]
  \| F_{4t}(t)\|^2 & \le C\delta_1 e^{-ct}.
 \end{array}
$$
For $F_{2t}$, 
$$
 \begin{array}{ll}
 \| F_{2t}(t)\|^2 & = \| (\alpha-\underline{\alpha})U_t\|^2 \le \| (\alpha-\underline{\alpha})V_{xxx}\|^2 \le C\delta_0(1+t)^{-5} \| x\cdot (\alpha-\underline{\alpha})\|^2 \\[5pt]
    &  \le C\delta_0 (1+t)^{-5} \ \ \mbox{if} \ \  x\cdot (\alpha-\underline{\alpha}) \in L^2,
 \end{array}
$$
because of $V_{xxx} = C\delta_0x\cdot ((1+t)^{-5/2}+ Cx^2 (1+t)^{-7/2}) e^{-\frac{x^2}{4\mu (1+t)}}$. Hence $\| F_{t}(t)\|^2 \le C\delta_1(1+t)^{-9/2}$.

The estimate of $\| F_{tx}(t)\|^2$ is done samely as above and omitted.
\\
\\
{\Large {\it Acknowledgment}.}  The authors would like to thank Professor Ming Mei who shortly visited Tokyo Institute of Technology in 2022-23. Through discussions with him they had a motive to consider the present problem.

\end{document}